\documentclass[12pt]{article}
\usepackage{amsmath,amsfonts,amssymb,amsthm}
\usepackage{color}
\usepackage{colortbl}
\usepackage{array}
\usepackage{mathrsfs}
\usepackage{dcolumn}
\usepackage{longtable}
\usepackage{hhline}
\usepackage{bm}
\newcommand{\rbeta}[2]{\genfrac{\{}{\}}{0pt}{0}{#1}{#2}}
\newcommand{\angled}[2]{\genfrac{\langle}{\rangle}{0pt}{0}{#1}{#2}}

\theoremstyle{plain}
\newtheorem{thm}{Theorem}[section]
\newtheorem{lem}[thm]{Lemma}
\newtheorem{cor}[thm]{Corollary}

\theoremstyle{definition}
\newtheorem{defn}[thm]{Definition}

\title{\bf Explicit Formula For Generalization Of Poly-Bernoulli Numbers and
Polynomials with $a,b,c$ Parameters}
\author{
{\large Hassan Jolany and Roberto B. Corcino}
}

\date{}

\begin{document}

\maketitle

\begin{abstract}
In this paper we investigate special generalized Bernoulli polynomials with $a, b, c$ parameters
that generalize classical Bernoulli numbers and polynomials. The present paper deals
with some recurrence formulae for the generalization of poly-Bernoulli numbers and polynomials
with $a, b, c$ parameters. Poly-Bernoulli numbers satisfy certain recurrence relationships
which are used in many computations involving poly-Bernoulli numbers. Obtaining a closed
formula for generalization of poly-Bernoulli numbers with $a, b, c$ parameters therefore seems
to be a natural and important problem. By using the generalization of poly-Bernoulli polynomials
with $a, b, c$ parameters of negative index we define symmetrized generalization of
poly-Bernoulli polynomials with a; b; c parameters of two variables and we prove duality
property for them. Also by Stirling numbers of the second kind we will find a closed formula
for them. Furthermore we generalize the Arakawa-Kaneko Zeta functions and by using the
Laplace-Mellin integral, define generalization of Arakawa-Kaneko Zeta functions with
$a, b, c$ parameters and obtain an interpolation formula for the generalization of poly-
Bernoulli numbers and polynomials with $a, b, c$ parameters. Furthermore we present a link
between this type of Zeta functions and Dirichlet series. By our interpolation formula, we
will interpolate the generalization of Arakawa-Kaneko Zeta functions with $a, b, c$ parameters.   


\end{abstract}

\section{Introduction}
The poly-Bernoulli polynomials have been studied by many researchers in recent decade.
The poly-Bernoulli polynomials have wide-ranging applications from number theory and combinatorics
to other fields of applied mathematics. One of applications of poly-Bernoulli
numbers that was investigated by Chad Brewbaker in \cite{Brew1, Brew2}, is about the number of (0; 1)-
matrices with n-rows and k columns. He showed the number of (0, 1)-matrices with n-rows
and k columns uniquely reconstructable from their row and column sums are the poly-
Bernoulli numbers of negative index $B_n^{(k)}$. Another application of poly-Bernoulli numbers
is in Zeta function theory. Multiple Zeta functions at non-positive integers can be described
in terms of these numbers. A third application of poly-Bernoulli numbers that was proposed
by Stephane Launois in \cite{Lau1, Lau2}, is about cardinality of some subsets of Sn. He proved
the cardinality of sub-poset of the reverse Bruhat ordering is equal to the poly-Bernoulli
numbers. Also one of other applications of poly-Bernoulli numbers is about skew Ferrers
boards. In \cite{Jona}, Jonas Sjostran found a relation between poly-Bernoulli numbers and the
number of elements in a Bruhat interval. Also he showed the Poincare polynomial (for value
q = 1) of some particularly interesting intervals in the finite Weyl group can be written in
terms of poly-Bernoulli numbers. Moreover Peter Cameron in \cite{Came} showed that the number
of acyclic orientations of a complete bipartite graph is a poly-Bernoulli number.

\smallskip
One of generalizations of poly-Bernoulli numbers that was first proposed by Y. Hamahata,
is the Multi-poly-Bernoulli numbers and he derived a closed formula for them. A. Bayad  ,
introduced a new generalization of poly-Bernoulli numbers and polynomials. He, by using
Dirichlet character, defined generalized poly-Bernoulli numbers associated to $\chi$. Also, he
introduced the generalized Arakawa-Kaneko $L$-functions and showed that the non-positive
integer values of the complex variable s of these $L$-functions can be written rationally in
terms of generalized poly-Bernoulli polynomials associated to $\chi$.

\smallskip
In \cite{Kim1, Kim2}, D. S. Kim and T. Kim considered poly-Bernoulli mixed-type polynomials.
From the properties of Sheffer sequences of these polynomials arising from umbrral calculus,
they derived several new and interesting identities. Also they introduced new generating
function which is known as Hermite poly-Bernoulli mixed-type polynomials.

\smallskip
In \cite{Jola}, H. Jolany et al, by using real $a, b, c$ parameters, introduced the generalization
of poly-Bernoulli polynomials with a; b; c parameters and found a closed relationships between
generalized poly-Bernoulli polynomials with a; b; c parameters and generalized Euler
polynomials with a; b; c parameters.

\smallskip
Let us briefly recall poly-Bernoulli numbers and polynomials. For an integer $k\in\mathbb{Z}$
\begin{equation}\label{eq1}
{\rm Li}_k(z)=\sum_{n=0}^{\infty}\frac{z^n}{n^k}
\end{equation}
which is the $k$-th polylogarithm if $k\geq 1$, and a rational function if $k\le0$. The name of the
function comes from the fact that it may alternatively be defined as the repeated integral of
itself, namely that
\begin{equation}\label{eq2}
{\rm Li}_{k+1}(z)=\int_{0}^{z}\frac{{\rm Li}_k(t)}{t}dt.
\end{equation}
One knows that ${\rm Li}_1 (z) = -\log (1-z)$. Also if $k$ is a negative integer, say $k = -r$, then
the poly-logarithmic function converges for $|x| < 1$ and equals
\begin{equation}\label{eq3}
{\rm Li}_{-r}(x)=\frac{\sum_{j=0}^{r}\angled{r}{j}x^{r-j}}{(1-x)^{r-j}}
\end{equation}
where the $\angled{r}{j}$ are the Eulerian numbers. The Eulerian numbers $\angled{r}{j}$ are the number of
permutations of $\{1, 2, \ldots, r\}$ with $j$ permutation ascents. One has
\begin{equation}\label{eq4}
\angled{r}{j}=\sum_{l=0}^{r+1}(-1)^{l}\binom{r+1}{l}(j-l+1)^r.
\end{equation}
The formal power series Lik (z) can be used to define poly-Bernoulli numbers and polynomials. The polynomials $B^{(k)}_n (x)$, 
$(n = 0, 1, 2, \ldots  )$ are said to be poly-Bernoulli polynomials if they satisfy
\begin{equation}\label{eq5}
\frac{{\rm Li}_k(1-e^{-t})}{1-e^{-t}}e^{xt}=\sum_{n=0}^{\infty}B^{(k)}_n(x)\frac{t^n}{n!}
\end{equation}
where $k\geq 1$. By (\ref{eq2}), the left-hand side of (\ref{eq5}) can be written in the form of iterated
integrals
\begin{equation}\label{eq6}
e^x\frac{1}{1+e^{x}}\int_{0}^{x}\frac{1}{1+e^{x}}\int_{0}^{x}\ldots\frac{1}{1+e^{x}}\int_{0}^{x}\frac{1}{1+e^{x}}dx\ldots dx=\sum_{n=0}^{\infty}B^{(k)}_n\frac{x^n}{n!}.
\end{equation}

For any $n\geq0$, we have
\[(-1)^n B^{(1)}_n(-x) = B_n(x)\]
where $Bn_ (x)$ are the classical Bernoulli polynomials given by
\begin{equation}\label{eq7}
\frac{t}{e^{t}-1}e^{xt}=\sum_{n=0}^{\infty}B_n(x)\frac{t^n}{n!}, \;\;|t|<2\pi.
\end{equation}
For $x = 0$ in (\ref{eq5}), we have $B^{(k)}_n (0) := B^{(k)}_n$, where $B^{(k)}_n$ are called poly-Bernoulli numbers
(for more information, see \cite{Kim2, San1, Adel, Copp, Kim3, Kama, Sasa, Shik, Ohno, Kim4}. In 2002, Q. M. Luo et al. \cite{Luo}, defined the generalization of Bernoulli numbers and polynomials with $a, b$ parameters as follows:
\begin{equation}\label{eq8}
\frac{t}{b^{t}-a^t}e^{xt}=\sum_{n=0}^{\infty}B_n(x;a,b)\frac{t^n}{n!}, \;\;\left|t\ln\frac{b}{a}\right|<2\pi.
\end{equation}
So, by (\ref{eq7}), we get
\[B_n(x; 1, e) := B_n(x), B_n(0; a, b) := B_n(a, b)\;\; \mbox{and}\;\; B_n(0;1, e) := B_n\]
where $B_n(a, b)$ are called the generalization of Bernoulli numbers with $a, b$ parameters. Also they the 
proved the following expression for this type of polynomials which interpolate the generalization of
Bernoulli polynomials with $a, b, c$ parameters
\[\sum_{j=1}^mj^n=\frac{1}{(n+1)(\ln b)^n}[B_{n+1} (m + 1; 1, b, b) - B_{n+1} (0; 1, b, b)].\]
H. Jolany et al. in \cite{Jola} defined a new generalization for poly-Bernoulli numbers and polynomials. They introduced the generalization of poly-Bernoulli polynomials with $a, b$ parameters as follows
\begin{equation}\label{dfn1}
\frac{{\rm Li}_k(1-(ab)^{-t})}{b^t-a^{-t}}e^{xt}=\sum_{n=0}^{\infty}B^{(k)}_n(x;a,b)\frac{t^n}{n!}.
\end{equation}
Also they extended the definition of generalized poly-Bernoulli polynomials with three parameters $a, b, c$ as follows:
\begin{equation}\label{dfn2}
\frac{{\rm Li}_k(1-(ab)^{-t})}{b^t-a^{-t}}c^{xt}=\sum_{n=0}^{\infty}B^{(k)}_n(x;a,b,c)\frac{t^n}{n!}.
\end{equation}
where $B^{(k)}_n (x; a, b, c)$ are called the generalization of poly-Bernoulli polynomials with $a, b, c$ parameters.
These are coefficients of power series expansion of a higher genus algabraic function with respect to a suitable variable. In the sequel, we list some closed formulas of poly-Bernoulli numbers and polynomials. 

\smallskip
Kim, in \cite{Kim1, Kim2, Kim3}, presented the following explicit formulas for poly-Bernoulli numbers 
\begin{align*}
B^{(k)}_n&=\frac{1}{n+1}\left\{B^{(k-1)}_n-\sum_{m=1}^{n-1}\binom{n}{m-1}B^{(k)}_m\right\}\\
B^{(-k)}_n&=\sum_{j=0}^{min(n,k)}(j!)^2\rbeta{n+1}{j+1}\rbeta{k+1}{j+1}, \;\;n, k\geq0,
\end{align*}
where
\[\rbeta{n}{m}=\frac{(-1)^m}{m!}\sum_{l=0}^m(-1)^l\binom{m}{l}l^n, \;\;n, m\geq0\]
called the Stirling numbers of the second kind. A. Bayad in \cite{Bay} introduced the generalized poly-Bernoulli polynomials $B^{(k)}_{n,\chi}(x)$. So, by applying their method, we introduce a closed formula and also interpolation formula for the generalization of poly-Bernoulli numbers and polynomials with a; b parameters which yields a deeper insight into the effectiveness of this type of generalizations.

\section{Explicit Formulas for Generalization of Poly-Bernoulli Polynomials with Three Parameters}

Now, we are in a position to state and prove the main results of this paper. In this section, we obtain some interesting new relations associated to generalization of poly-Bernoulli numbers and polynomials with $a, b, c$ parameters. Here we prove a collection of important
and fundamental identities involving this type of number and polynomials. We also deduce their special cases which leads to the corresponding results for the poly-Bernoulli polynomials. 

\smallskip
First of all, we present an explicit formula for generalization of poly-Bernoulli polynomials with $a, b, c$ parameters

\bigskip
\begin{thm}\label{th1} {\rm ({\bf Explicit Formula})} For $k \in \mathbb{Z}$, $n\geq0$, we have
\begin{equation}
B^{(k)}_n (x; a, b) =\sum_{m=0}^n\frac{1}{(m + 1)^k}\sum_{j=0}^m(-1)^j\binom{m}{j}(x - j \ln a - (j + 1) \ln b)^n.
\end{equation}
\begin{proof}
\begin{align*}
\frac{{\rm Li}_k(1-(ab)^{-t})}{b^t-a^{-t}}&=b^{-t}\left(\sum_{m=1}^{\infty}\frac{(1-(ab)^{-t})^{m-1}}{m^k}\right)=b^{-t}\left(\sum_{m=0}^{\infty}\frac{(1-(ab)^{-t})^{m}}{(m+1)^k}\right)\\
&=b^{-t}\sum_{m=0}^{\infty}\frac{1}{(m+1)^k}\sum_{j=0}^m(-1)^j\binom{m}{j}e^{-jt\ln (ab)}\\
&=\sum_{m=0}^{\infty}\frac{(1}{(m+1)^k}\sum_{j=0}^m(-1)^j\binom{m}{j}e^{-t(j\ln a+(j+1)\ln b)}
\end{align*}
So, we get
\begin{align*}
\frac{{\rm Li}_k(1-(ab)^{-t})}{b^t-a^{-t}}e^{xt}&=\sum_{m=0}^{\infty}\frac{1}{(m+1)^k}\sum_{j=0}^m(-1)^j\binom{m}{j}e^{t(x-j\ln a-(j+1)\ln b)}\\
&=\sum_{n=0}^{\infty}\left(\sum_{m=0}^{\infty}\frac{1}{(m+1)^k}\sum_{j=0}^m(-1)^j\binom{m}{j}(x-j\ln a-(j+1)\ln b)^n\right)\frac{t^n}{n!}
\end{align*}
By comparing the coefficients of $\frac{t^n}{n!}$ on both sides, the proof is completed.
\end{proof}
\end{thm}

As a direct result, by applying the same method as Theorem \ref{th1}, we derive following corollaries.

\bigskip
\begin{cor}\label{cor1} For $k \in \mathbb{Z}$, $n\geq0$, we have
\begin{equation}
B^{(k)}_n (x; a, b,c) =\sum_{m=0}^n\frac{1}{(m + 1)^k}\sum_{j=0}^m(-1)^j\binom{m}{j}(x\ln c - j \ln a - (j + 1) \ln b)^n.
\end{equation}
\end{cor}

As a direct result, by applying $a = e, b = 1, c = e$ in Corollary 1, we get the following corollary.

\bigskip
\begin{cor}\label{cor2} For $k \in \mathbb{Z}$, $n\geq0$, we have
\begin{equation}
B^{(k)}_n (x) =\sum_{m=0}^n\frac{1}{(m + 1)^k}\sum_{j=0}^m(-1)^j\binom{m}{j}(x - j)^n.
\end{equation}
\end{cor}

Furthermore, by setting k = 1 in Corollary 2 and because we have, $B_n (x) = (-1)^n B^{(1)}_n (-x)$, we obtain following explicit formulas for classical Bernoulli numbers and polynomials.

\bigskip
\begin{cor}\label{cor3} For $k \in \mathbb{Z}$, $n\geq0$, we have
\begin{align*}
B_n (x) &=\sum_{m=0}^n\frac{1}{m + 1}\sum_{j=0}^m(-1)^j\binom{m}{j}(x + j)^n\\
B_n&=\sum_{m=0}^n\frac{1}{m + 1}\sum_{j=0}^m(-1)^j\binom{m}{j}j^n.
\end{align*}
\end{cor}

Now, we investigate some recursive formulas for the generalization of poly-Bernoulli numbers and polynomials with $a, b$ parameters.

\bigskip
\begin{thm}\label{th2} {\rm ({\bf Recursive Formula})} For all $k\geq1$ and $n\geq0$, we have
\begin{equation}
B^{(k)}_n (x; a, b) =(\ln a + \ln b)\sum_{m=0}^n(-ln a)^m\binom{n}{m}B^{(k-1)}_{n-m} (a, b)\sum_{l=0}^m\frac{(-ln a)^{-l}}{n-l+1}
\binom{m}{l}B_{l}(x; a^{-1}, b).
\end{equation}
\begin{proof} We know
\[{\rm Li}_{k+1} (t) =\int_0^t\frac{{\rm Li}_k (s)}{s}ds\]
so
\[{\rm Li}_{k+1} (1-(ab)^{-t}) =\int_0^t\frac{{\rm Li}_k (1-(ab)^{-s})}{1-(ab)^{-s}}(ln ab) e^{-s \ln ab} ds.\]
So we get
\[\frac{{\rm Li}_{k+1} (1-(ab)^{-t})}{b^t-a^{-t}}e^{xt} =\frac{a^te^{xt}}{(ab)^t-1}\int_0^t(ln ab)\frac{{\rm Li}_k (1-(ab)^{-s})}{1-(ab)^{-s}}e^{-s \ln ab} ds.\]
Therefore, we obtain
\[\sum_{n=0}^{\infty}B^{(k)}_n(x;a,b)\frac{t^n}{n!}\qquad\qquad\qquad\qquad\qquad\qquad\qquad\qquad\qquad\qquad\qquad\qquad\qquad\qquad\qquad\qquad\qquad\qquad\qquad\qquad\]
\begin{align*}
&=(\ln ab)\left(\sum_{n=0}^{\infty}B^{(k)}_n(x;a^{-1},b)\frac{t^{n-1}}{n!}\right)\int_0^t\frac{{\rm Li}_k (1-(ab)^{-s})}{1-(ab)^{-s}}e^{-s \ln ab} ds\\
&=(\ln ab)\left(\sum_{n=0}^{\infty}B^{(k)}_n(x;a^{-1},b)\frac{t^{n-1}}{n!}\right)\int_0^t\left(\sum_{n=0}^{\infty}\frac{(-s\ln a)^n}{n!}\right)\left(\sum_{n=0}^{\infty}B^{(k-1)}_n(a,b)\frac{s^n}{n!}\right)ds\\
&=(\ln ab)\left(\sum_{n=0}^{\infty}B^{(k)}_n(x;a^{-1},b)\frac{t^{n-1}}{n!}\right)\sum_{n=0}^{\infty}\left(\sum_{m=0}^{n}(-\ln a)^{n-m}\binom{n}{m}B^{(k-1)}_m(a,b)\right)\frac{t^{n+1}}{(n+1)!}\\
&=(\ln ab)\sum_{n=0}^{\infty}\left(\sum_{l=0}^nB^{(k)}_{n-l}(x;a^{-1},b)\sum_{m=0}^{l}(-\ln a)^{l-m}\binom{l}{m}B^{(k-1)}_m(a,b)\frac{t^n}{(l+1)!(n-l)!}\right)\\
&=(\ln ab)\sum_{n=0}^{\infty}\left(\sum_{l=0}^n\frac{B^{(k)}_{n-l}(x;a^{-1},b)}{l+1}\binom{n}{l}\sum_{m=0}^{l}(-\ln a)^{l-m}\binom{l}{m}B^{(k-1)}_m(a,b)\right)\frac{t^n}{n!}
\end{align*}
So, by applying the following identity
\[\binom{n}{l}\binom{l}{m}=\binom{n}{m}\binom{n-m}{n-l},\]
we obtain
\[\sum_{n=0}^{\infty}B^{(k)}_n(x;a,b)\frac{t^n}{n!}\qquad\qquad\qquad\qquad\qquad\qquad\qquad\qquad\qquad\qquad\qquad\qquad\qquad\qquad\qquad\qquad\qquad\qquad\qquad\qquad\]
\begin{equation*}
=(\ln ab)\sum_{n=0}^{\infty}\left(\sum_{m=0}^nB^{(k-1)}_m(a,b)\binom{n}{m}\sum_{l=m}^{n}\frac{(-\ln a)^{l-m}}{l+1}\binom{n-m}{n-l}{B^{(k)}_{n-l}(x;a^{-1},b)}\right)\frac{t^n}{n!}.
\end{equation*}
Putting $l'=n-l$, we have
\[\sum_{n=0}^{\infty}B^{(k)}_n(x;a,b)\frac{t^n}{n!}\qquad\qquad\qquad\qquad\qquad\qquad\qquad\qquad\qquad\qquad\qquad\qquad\qquad\qquad\qquad\qquad\qquad\qquad\qquad\qquad\]
\begin{equation*}
=(\ln ab)\sum_{n=0}^{\infty}\left(\sum_{m=0}^nB^{(k-1)}_m(a,b)\binom{n}{m}\sum_{l'=0}^{n-m}\frac{(-\ln a)^{n-l'-m}}{n-l'+1}\binom{n-m}{l'}{B^{(k)}_{l'}(x;a^{-1},b)}\right)\frac{t^n}{n!}.
\end{equation*}
Putting $m'=n-m$, we obtain
\[\sum_{n=0}^{\infty}B^{(k)}_n(x;a,b)\frac{t^n}{n!}\qquad\qquad\qquad\qquad\qquad\qquad\qquad\qquad\qquad\qquad\qquad\qquad\qquad\qquad\qquad\qquad\qquad\qquad\qquad\qquad\]
\begin{equation*}
=(\ln ab)\sum_{n=0}^{\infty}\left(\sum_{m'=0}^nB^{(k-1)}_{n-m'}(a,b)\binom{n}{m'}\sum_{l'=0}^{m'}\frac{(-\ln a)^{m'-l'}}{n-l'+1}\binom{m'}{l'}{B^{(k)}_{l'}(x;a^{-1},b)}\right)\frac{t^n}{n!}.
\end{equation*}
By comparing the coefficients of $\frac{t^n}{n!}$ on both sides, the proof is complete.
\end{proof}
\end{thm}

As a direct consequence of Theorem \ref{th2} with $a = e, b = 1$, we obtain the following corollary
which is the well known recurrence formula for classical poly-Bernoulli polynomials.

\bigskip
\begin{cor}\label{cor4}
For all $k \geq 1$, $n \geq 0$, we have
\begin{equation}
B^{(k)}_n (x) =\sum_{m=0}^n(-1)^m\binom{n}{m}B^{(k-1)}_{n-m}\sum_{l=0}^m\frac{(-1)^{l}}{n-l+1}
\binom{m}{l}B^{(1)}_{l}(x).
\end{equation}
\end{cor}

Let us consider the extreme recurrence formula for generalization of poly-Bernoulli polynomials
with $a, b$ parameters. By using following lemma and some standard techniques based
upon generating function and series rearrangement we present a new recurrence formula for
generalization of poly-Bernoulli polynomials with $a, b$ parameters.

\bigskip
\begin{lem}\label{lem1}
For $a, b > 0$ and $n \geq0$, we have
\begin{equation}
B^{(k)}_n (x; a, b) =(\ln a + \ln b)^nB^{(k)}_n \left(\frac{x-\ln b}{\ln a + \ln b}\right).
\end{equation}
\begin{proof}
By applying (\ref{dfn1}), we have
\begin{align*}
\sum_{n=0}^{\infty}B^{(k)}_n(x;a,b)\frac{t^n}{n!}&=\frac{{\rm Li}_k(1-(ab)^{-t})}{b^t-a^{-t}}e^{xt}=\frac{1}{b^t}\frac{{\rm Li}_k(1-(ab)^{-t})}{1-(ab)^{-t}}e^{xt}\\
&=e^{(x-\ln b)t}\frac{{\rm Li}_k(1-e^{-t\ln ab}}{1-e^{-t\ln ab}}=\sum_{n=0}^{\infty}(\ln a + \ln b)^nB^{(k)}_n \left(\frac{x-\ln b}{\ln a + \ln b}\right)\frac{t^n}{n!}
\end{align*}
So, by comparing the coefficients of $\frac{t^n}{n!}$ on both sides, we obtain the desired result.
\end{proof}
\end{lem}

Now we are ready to present our second recurrence formula for generalization of Poly-Bernoulli numbers and polynomials with $a, b$ parameters.

\bigskip
\begin{thm}\label{th3}
For $k \in \mathbb{Z}$ and $n \geq 2$, we have
\begin{align*}
B^{(k)}_0(x;a,b)&=1\\
B^{(k)}_1(x;a,b)&=\frac{1}{2}\left[B^{(k-1)}_1(x; a, b) + \left(\frac{x-\ln b}{\ln a + \ln b}\right)B^{(k)}_0(x;a,b)\right]\\
B^{(k)}_n (x; a, b)&=\frac{1}{n+1}\left\{B^{(k-1)}_n (x; a, b) + (x - \ln b) (ln a + ln b)^{n-1} B^{(k)}_0 (x; a, b)\right.\\
&\;\;\;\;+(x - \ln b)\sum_{m=1}^{n-1}(\ln a + \ln b)^{n-m-1}\binom{n}{m}B^{(k)}_m (x; a, b)\\
&\;\;\;\;-\left.\sum_{m=1}^{n-1}(\ln a + \ln b)^{n-m}\binom{n}{m-1}B^{(k)}_m (x; a, b)\right\}
\end{align*}
\begin{proof}
From [], we have the following recurrence formula for poly-Bernoulli polynomials
\begin{equation}\label{eqn5}
B^{(k)}_n (x; a, b)=\frac{1}{n+1}\left[B^{(k-1)}_n (x) + xB^{(k)}_0 (x)\sum_{m=1}^{n-1}\left[\binom{n}{m-1}-\binom{n}{m}x\right]B^{(k)}_m (x)\right]
\end{equation}
So, by applying Lemma \ref{lem1} and replacing $r$ by $\frac{x-\ln b}{\ln a + \ln b}$ in (\ref{eqn5}), we obtain the desired result.
\end{proof}
\end{thm}

Now, we show that the generalization of poly-Bernoulli polynomials of a; b parameters are in the set of Appell polynomials.

\smallskip
For a sequence $\{P_n (x)\}_{n=0}^{\infty}$ of Appell polynomials, which is a sequence of polynomials satisfying
\[\frac{dP_n (x)}{dx}=nP_{n-1}(x), \;\;n\geq1.\]
Tremendous properties are well known. Among them, the most important classifications of Appell polynomials may be the following equivalent conditions (\cite{Lee, Shoh, Tosc}).

\bigskip
\begin{thm}\label{th5}
Let $\{P_n (x)\}_{n=0}^{\infty}$ be a sequence of polynomials. Then the following are all equivalent

\smallskip
(a) $\{P_n (x)\}_{n=0}^{\infty}$ is a sequence of Appell polynomials.

\smallskip
(b) $\{P_n (x)\}_{n=0}^{\infty}$ has a generating function of the form
\[A(t)e^{xt}=\sum_{n=0}^{\infty}P_n(x)\frac{t^n}{n!},\]

\smallskip
where $A(t)$ is a formal power series in $t$ with $A(0) \neq 0$.

\smallskip
(c) $\{P_n (x)\}_{n=0}^{\infty}$ satisfies 
\[ P_n (x + y) =\sum_{k=0}^n\binom{n}{k}P_{n-k} (x) y^k\]
\end{thm}

Now, in the following theorem we prove that the generalization of poly-Bernoulli polynomials are in the set of Appell sequence

\begin{thm}\label{th6}
{\rm ({\bf Appell Sequence})} The generalized poly-Bernoulli polynomials satisfy the following differential equation
\begin{align}
\frac{dB^{(k)}_0 (x; a, b)}{dx}&=0\nonumber\\
\frac{dB^{(k)}_{n+1} (x; a, b)}{dx}&=(n + 1)B^{(k)}_n (x; a, b).
\end{align}
\begin{proof}
By differentiating both sides of (\ref{dfn1}), with respect to $x$, we have
\[t\frac{{\rm Li}_{k}(1-(ab)^{-t})}{b^t-a^{-t}}e^{xt}=\sum_{n=0}^{\infty}\frac{dB^{(k)}_{n} (x; a, b)}{dx}\frac{t^n}{n!}\]
and obtain
\[\frac{{\rm Li}_{k}(1-(ab)^{-t})}{b^t-a^{-t}}e^{xt}=\sum_{n=0}^{\infty}\left[\frac{1}{(n+1)}\frac{dB^{(k)}_{n+1} (x; a, b)}{dx}\right]\frac{t^n}{n!}\]
which yields the desired results.
\end{proof}
\end{thm}
Thus, by applying the property of (c) of Theorem \ref{th5}, we obtain following corollary.

\begin{cor}
{\rm ({\bf Addition Formula})} For $k \in \mathbb{Z}$ and $n \geq 0$, we have
\begin{equation}
B^{(k)}_n (x + y; a, b) =\sum_{m=0}^n\binom{n}{m}B^{(k)}_m (x; a, b) y^{n-m}.
\end{equation}
\end{cor}

\bigskip
\noindent In particular,
\begin{equation}
B^{(k)}_n (x; a, b) =\sum_{m=0}^n\binom{n}{m}B^{(k)}_m (a, b) x^{n-m}.
\end{equation}
and by taking $y = (m - 1) x$, we obtain Multiplication theorem for them
\begin{equation}
B^{(k)}_n (mx; a, b) =\sum_{i=0}^n\binom{n}{i}B^{(k)}_i (x; a, b) (m-1)^{n-i}x^{n-i},\;\;m=1, 2, \ldots
\end{equation}

\smallskip
Actually, because generalization of poly-Bernoulli polynomials of $a, b$ parameters are in
the set of Appell polynomials, we can derive numerous properties for them. For instance in
\cite{Cost}, F. A. Costabile and E. Longo presented a new definition by means of a determinantal
form for Appell polynomials by using of linear algebra tools and also M. E. H. Ismail in \cite{Isma},
found a differential equation for Appell polynomials.

\section{Symmetrized Generalization of Poly-Bernoulli Polynomials with $a, b$ parameters of Two Variables}

Kaneko, Japanese mathematician introduced the symmetrized poly-Bernoulli polynomials with two variables and by using their method we introduce symmetrized generalization of poly-Bernoulli polynomials with $a, b$ parameters of two variables and construct a generating
function for symmetrized generalization of poly-Bernoulli polynomials with $a, b$ parameters of two variables. Also we give a closed formula and duality property for this type of polynomials as well.

\begin{defn}
For $m, n \geq 0$, we define
\begin{equation}
C^{(-m)}_n (x, y; a, b) =\frac{1}{(\ln a + \ln b)^n}\sum_{k=0}^m\binom{m}{k}B^{(-k)}_n (x; a, b)\left(y-\frac{\ln b}{\ln a + \ln b}\right)^{m-k}
\end{equation}
\end{defn}

Now, in the following theorem we introduce a generating function for $C^{(-m)}_n (x, y; a, b)$.

\begin{thm}\label{th7}
For $m, n \geq 0$, we have
\begin{equation}
\sum_{n=0}^{\infty}\sum_{m=0}^{\infty}C^{(-m)}_n (x, y; a, b)\frac{t^n}{n!}\frac{u^m}{m!}=\frac{e^{\left(x+\frac{\ln a}{\ln a+\ln b}\right)t}e^{\left(y+\frac{\ln b}{\ln a+\ln b}\right)u}}{e^t + e^u - e^{t+u}}
\end{equation}
\begin{proof}
By using the definition of $C^{(-m)}_n (x, y; a, b)$, the left-hand side can be written as
\[LHS=\sum_{n=0}^{\infty}\sum_{m=0}^{\infty}\frac{1}{(\ln a+\ln b)^n}\sum_{k=}^mB^{(-k)}_n (x; a, b)\left(y-\frac{\ln b}{\ln a+\ln b}\right)^{m-k}\frac{t^n}{n!}\frac{u^m}{k!(m-k)!}\]
By putting $l = m - k$, we get
\begin{align*}
LHS &= \sum_{n=0}^{\infty}\sum_{k=0}^{\infty}\sum_{l=0}^{\infty}\frac{1}{(\ln a + \ln b)^n}B^{(-k)}_n (x; a, b)\left(y-\frac{\ln b}{\ln a+\ln b}\right)^{l}\frac{t^n}{n!}\frac{u^k}{k!}\frac{u^l}{l!}\\
&= e^{\left(y-\frac{\ln b}{\ln a+\ln b}\right)u}\sum_{n=0}^{\infty}\sum_{k=0}^{\infty}\frac{1}{(\ln a + \ln b)^n}B^{(-k)}_n (x; a, b)\frac{t^n}{n!}\frac{u^k}{k!}\\
&= e^{\left(y-\frac{\ln b}{\ln a+\ln b}\right)u}\sum_{k=0}^{\infty}\left(e^{xt}\sum_{n=0}^{\infty}B^{(-k)}_n (a, b)\frac{\left(\frac{t}{\ln a+\ln b}\right)^n}{n!}\right)\frac{u^k}{k!}\\
&= e^{\left(y-\frac{\ln b}{\ln a+\ln b}\right)u}\sum_{k=0}^{\infty}\left(e^{xt}\frac{{\rm Li}_{-k} (1 - e^{-t})}{1-e^{-t}}e^{\left(\frac{-t\ln b}{\ln a+ \ln b}\right)^n}\right)\frac{u^k}{k!}\\
&= e^{\left(y-\frac{\ln b}{\ln a+\ln b}\right)u}e^{\left(x-\frac{\ln b}{\ln a+ \ln b}\right)t}\sum_{k=0}^{\infty}\sum_{n=0}^{\infty}B^{(-k)}_n\frac{t^n}{n!}\frac{u^k}{k!}
\end{align*}
But Kaneko proved following expression
\[\sum_{k=0}^{\infty}\sum_{n=0}^{\infty}B^{(-k)}_n\frac{t^n}{n!}\frac{u^k}{k!}=\frac{e^{t+u}}{e^t + e^u - e^{t+u}}.\]
So, by applying this expression, we obtain the desired result.
\end{proof}
\end{thm}

\smallskip
As a direct result, we have the following corollary for $C^{(-m)}_n (x, y; a, b)$ that is the well known duality property.

\bigskip
\begin{cor}
{\rm ({\bf Duality Property})} For $m \geq 0$, we have
\begin{equation}
C^{(-m)}_n (x, y; a, b)=C^{(-m)}_n (y, x; b, a).
\end{equation}
\end{cor}

Now, we are ready to show a closed formula for $C^{(-m)}_n (x, y; a, b)$ which is important and fundamental.

\bigskip
\begin{thm}\label{th8}
{\rm ({\bf Closed Formula})} For $m \geq 0$, we have
\begin{align}
C^{(-m)}_n (x, y; a, b)&=\sum_{j=0}^{\infty}(j!)^2\left(\sum_{p=0}^{\infty}\left(x+\frac{\ln a}{\ln a+ \ln b}\right)^{n-p}\binom{n}{p}\rbeta{p}{j}\right)\times\\
&\;\;\;\times\left(\sum_{l=0}^{\infty}\left(y+\frac{\ln b}{\ln a+ \ln b}\right)^{m-l}\binom{m}{l}\rbeta{l}{j}\right)\nonumber
\end{align}
\begin{proof}
By applying Theorem \ref{th7}, we have
\begin{align*}
\sum_{n=0}^{\infty}\sum_{m=0}^{\infty}C^{(-m)}_n (x, y; a, b)\frac{t^n}{n!}\frac{u^m}{m!}&=\frac{e^{\left(x+\frac{\ln a}{\ln a+\ln b}\right)t}e^{\left(y+\frac{\ln b}{\ln a+\ln b}\right)u}}{e^t + e^u - e^{t+u}}=\frac{e^{\left(x+\frac{\ln a}{\ln a+\ln b}\right)t}e^{\left(y+\frac{\ln b}{\ln a+\ln b}\right)u}}{1 - (e^t - 1) (e^u - 1)}\\
&=e^{\left(x+\frac{\ln a}{\ln a+\ln b}\right)t}e^{\left(y+\frac{\ln b}{\ln a+\ln b}\right)u}\sum_{j=0}^{\infty}(e^t - 1)^j (e^u - 1)^j\\
&=\sum_{j=0}^{\infty}e^{\left(x+\frac{\ln a}{\ln a+\ln b}\right)t}(e^t - 1)^j e^{\left(y+\frac{\ln b}{\ln a+\ln b}\right)u}(e^u - 1)^j
\end{align*}
By applying the generating function of Stirling numbers of second kind
\[\sum_{n=0}^{\infty}\rbeta{n}{k}\frac{u^n}{n!}=\frac{(e^u-1)^k}{k!}\]
the right-hand side of the last expression becomes
\begin{align*}
&=\sum_{j=0}^{\infty}\left(j!\sum_{n=0}^{\infty}\frac{\left(x+\frac{\ln a}{\ln a+\ln b}\right)^nt^n}{n!}\sum_{m=0}^{\infty}\rbeta{m}{j}\frac{t^m}{m!}\right)\left(j!\sum_{n=0}^{\infty}\frac{\left(y+\frac{\ln b}{\ln a+\ln b}\right)^nu^n}{n!}\sum_{m=0}^{\infty}\rbeta{m}{j}\frac{u^m}{m!}\right)\\
&=\sum_{j=0}^{\infty}\left(j!\sum_{l=0}^{\infty}\sum_{m=0}^{l}\left(x+\frac{\ln a}{\ln a+\ln b}\right)^{l-m}\binom{l}{m}\rbeta{m}{j}\frac{t^l}{l!}\right)\times\\
&\;\;\;\times\left(j!\sum_{p=0}^{\infty}\sum_{r=0}^{p}\left(y+\frac{\ln b}{\ln a+\ln b}\right)^{p-r}\binom{p}{r}\rbeta{r}{j}\frac{u^p}{p!}\right)\\
&=\sum_{l=0}^{\infty}\sum_{p=0}^{\infty}\frac{t^l}{l!}\frac{u^p}{p!}\sum_{j=0}^{\infty}(j!)^2\left(\sum_{m=0}^{l}\left(x+\frac{\ln a}{\ln a+\ln b}\right)^{l-m}\binom{l}{m}\rbeta{m}{j}\right)\times\\
&\;\;\;\times\left(\sum_{r=0}^{p}\left(y+\frac{\ln b}{\ln a+\ln b}\right)^{p-r}\binom{p}{r}\rbeta{r}{j}\right)
\end{align*}
which yields the result.
\end{proof}
\end{thm}

\section{Generalization of Arakawa-Kaneko $L$-Functions with $a, b$ Parameters}

It is well known since the second-half of the 19-th century the Riemann Zeta function may be represented by the normalized Mellin transformation
\[\zeta(s)=\frac{1}{\Gamma(s)}\int_{0}^{\infty}t^{s-1}\frac{e^{-t}}{1-e^{-t}}dt, \;\;{\rm Re}(s)>1.\]
T. Arakawa and M. Kaneko, by inspiration of last expression, introduced Arakawa-Kaneko Zeta function as follows. For any integer $k \geq 1$
\[\xi_k(s, x)=\frac{1}{\Gamma(s)}\int_{0}^{\infty}\frac{{\rm Li}_k(1-e^{-t})}{1-e^{-t}}e^{-xt}t^{s-1}dt.\]
It is defined for ${\rm Re} (s) > 0$ and $x > 0$ if $k \geq 1$, and for ${\rm Re} (s) > 0$ and $x > |k|+1$ if $k \le 0$.
The function $\xi_k (s, x)$ has analytic continuation to an entire function on the whole complex $s$-plane and
\begin{equation}\label{eq41}
\xi_k(-n, x) = (-1)^n B^{(k)}_n (-x)
\end{equation}
for all non-negative integer $n$ and $x \geq 0$ (for more information, see [7]). 

\smallskip
For $k \in \mathbb{Z}$, the generalization of Arakawa-Kaneko Zeta function with $a, b$ parameters are given by the Laplace Mellin-integral
\begin{equation}
\xi_k(s, x; a, b)=\frac{1}{\Gamma(s)}\int_{0}^{\infty}\frac{{\rm Li}_k(1-(ab)^{-t})}{b^t-a^{-t}}e^{-xt}t^{s-1}dt.
\end{equation}
It is defined for ${\rm Re} (s) > 0$ and $x > 0$ if $k \geq 1$, and for ${\rm Re} (s) > 0$ and $x > |k| + 1$ if $k \le 0$. It is easy to see that the generalization of Arakawa-Kaneko zeta function with $a, b$ parameters include the Arakawa-Kaneko Zeta function and Hurwitz-Zeta function.

\smallskip
In this section, we now derive an interpolation formula of generalization of poly-Bernoulli polynomials with $a, b$ parameters and investigate fundamental properties of $\xi_k (s, x; a, b)$. At first, in following lemma we give a relation between generalization of Arakawa-Kaneko Zeta function with $a, b$ parameters and classical Arakawa-Kaneko Zeta function.

\bigskip
\begin{lem}\label{lem2}
For $k \in \mathbb{Z}$, we have
\begin{equation}
\xi_k(s, x; a, b)=\frac{1}{(\ln a+\ln b)^s}\xi_k\left(s, \frac{x+\ln b}{\ln a+\ln b}\right). 
\end{equation}
\begin{proof}
It is easy to see that
\begin{equation*}
\xi_k(s, x; a, b)=\frac{1}{\Gamma(s)}\int_{0}^{\infty}\frac{{\rm Li}_k(1-e^{-t\ln ab})}{1-e^{-t\ln ab}}e^{-(x+\ln b)t}t^{s-1}dt.
\end{equation*}
So, by changing $t$ by $z = (\ln a + \ln b) t$, we obtain
\begin{equation*}
\xi_k(s, x; a, b)=\frac{1}{(\ln a + \ln b)^s}\frac{1}{\Gamma(s)}\int_{0}^{\infty}\frac{{\rm Li}_k(1-e^{-z})}{1-e^{-z}}e^{-\left(\frac{(x+\ln b)}{\ln a+ \ln b}\right)z}t^{s-1}dt,
\end{equation*}
which yields the lemma.
\end{proof}
\end{lem}

\bigskip
\begin{thm}\label{th9}
{\rm ({\bf Interpolation formula})} The function $s \to \xi_k (s, x; a, b)$ has analytic continuation to an entire function on the whole complex $s$-plane and for any positive integer $n$, we have
\begin{equation}
\xi_k(-n, x; a, b)=(-1)^n B^{(k)}_n (-x; a, b).
\end{equation}
\begin{proof}
By using Lemma \ref{lem2}, to prove that $s \to \xi_k (s, x; a, b)$ has analytic continuation to an entire function on the whole complex $s$-plane, it is sufficient to show that $s \to \xi_k (s, x)$ has such a property. Since this fact comes from the first part of Theorem 1.10 in [7] we omit it. By using Lemma \ref{lem1} and expression (\ref{eq41}), we get
\begin{align*}
\xi_k(-n, x; a, b)&=(\ln a + \ln b)^n\xi_k\left(-n,\frac{x + \ln b}{\ln a + \ln b}\right)\\
&=(-1)^n (\ln a + \ln b)^n B^{(k)}_n\left(\frac{-x - \ln b}{\ln a + \ln b}\right)\\
&=(-1)^n B^{(k)}_n (-x; a, b).
\end{align*}
So we obtain the desired result.
\end{proof}
\end{thm}

As an immediate consequence of previous theorems in this section, we obtain an explicit formula for $\xi_k (s, x; a, b)$.

\bigskip
\begin{cor}
For $k \in \mathbb{Z}$, we have
\begin{equation}
\xi_k (s, x; a, b)=\sum_{n=0}^{\infty}\frac{1}{(n+1)^k}\sum_{j=0}^n(-1)^j\binom{n}{j}\frac{1}{(x + j \ln a + (j + 1) \ln b)^s}.
\end{equation}
\begin{proof}
By applying Theorem \ref{th9}, we can interchange the integral and the sum. Hence
\begin{align*}
\xi_k(s, x; a, b)&=\frac{1}{\Gamma(s)}\sum_{n=1}^{\infty}\frac{1}{n^k}\int_{0}^{\infty}(1 - e^{-t \ln ab})^{n-1}e^{-t(x+ \ln b)}t^{s-1}dt\\
&=\frac{1}{\Gamma(s)}\sum_{n=0}^{\infty}\frac{1}{(n+1)^k}\int_{0}^{\infty}\sum_{j=0}^n\binom{n}{j}(-1)^je^{-t(x+j\ln a+(j+1) \ln b)}t^{s-1}dt\\
&=\sum_{n=0}^{\infty}\frac{1}{(n+1)^k}\sum_{j=0}^n\binom{n}{j}(-1)^j\frac{1}{\Gamma(s)}\int_{0}^{\infty}e^{-t(x+j\ln a+(j+1) \ln b)}t^{s-1}dt\\
&=\sum_{n=0}^{\infty}\frac{1}{(n+1)^k}\sum_{j=0}^n\binom{n}{j}(-1)^j\frac{1}{(x + j \ln a + (j + 1) \ln b)^s}.
\end{align*}
So the proof is complete.
\end{proof}
\end{cor}

Raabe's formula is a fundamental and universal property in the theory of Zeta function and plays an important role in special functions. Raabe's formula holds for several types of Zeta functions. For instance, Hurwitz Zeta function, Euler Zeta function and $q$-Euler Zeta function, multiple Zeta function. This formula provides a powerful link between zeta integrals and Dirichlet series. Raabe's formula can be obtained from the Hurwitz zeta function
\[\zeta(s,q)=\sum_{n=0}^{\infty}\frac{1}{(n+q)^s}\]
via the integral formula
\[\int_0^1\zeta(s, q + t) dq=\frac{t^{1-s}}{s - 1}.\]

\smallskip
Now, in next theorem, we will present a interesting link between integral of generalization of Arakawa-Kaneko zeta function with $a, b$ parameters and Dirichlet series. In fact we prove the Raabe's formula for our new types of zeta function.

\begin{lem}\label{lem3}
{\rm ({\bf Difference Formula})} we have
\[\xi_k (s, x + \ln ab; a, b) - \xi_k (s, x; a, b) \qquad\qquad\qquad\qquad\qquad\qquad\qquad\qquad\qquad\qquad\qquad\qquad\qquad\qquad\qquad\qquad\qquad\qquad\qquad\qquad\qquad\qquad\]
\begin{equation}
=\sum_{m=0}^{\infty}\frac{1}{(m+1)^k}\sum_{j=0}^{m+1}\binom{m+1}{j}(-1)^{j+1}(x + j \ln a + (j + 1) \ln b)^{-s}.
\end{equation}
\begin{proof}
By applying the definition of generalization of Arakawa-Kaneko zeta function with $a, b$ parameters, we get
\[\xi_k (s, x + \ln ab; a, b) - \xi_k (s, x; a, b) \qquad\qquad\qquad\qquad\qquad\qquad\qquad\qquad\qquad\qquad\qquad\qquad\qquad\qquad\qquad\qquad\qquad\qquad\qquad\qquad\qquad\qquad\]
\begin{align*}
&=-\frac{1}{\Gamma(s)}\int_0^{\infty}{\rm Li}_k(1 - e^{-t \ln ab})e^{-t(x+\ln b)}t^{s-1}dt\\
&=-\frac{1}{\Gamma(s)}\sum_{m=1}^{\infty}\frac{1}{m^k}\int_0^{\infty}(1 - (ab)^{-t})^me^{-t(x+\ln b)}t^{s-1}dt\\
&=\frac{1}{\Gamma(s)}\sum_{m=0}^{\infty}\frac{1}{(m+1)^k}\int_0^{\infty}\sum_{j=0}^{m+1}\binom{m+1}{j}(-1)^{j+1}e^{(x + j \ln a + (j + 1) \ln b)}t^{s-1}dt\\
&=\sum_{m=0}^{\infty}\frac{1}{(m+1)^k}\sum_{j=0}^{m+1}\binom{m+1}{j}(-1)^{j+1}(x + j \ln a + (j + 1) \ln b)^{-s}.
\end{align*}
So, we obtain the desired result.
\end{proof}
\end{lem}

Now, we are ready to present the Raabe's formula for our new types of zeta function.

\begin{thm}\label{th10}
{\rm ({\bf (Raabe's Formula})} we have
\[\int_{0}^{\ln ab}\xi_k (s, x + w; a, b) dw \qquad\qquad\qquad\qquad\qquad\qquad\qquad\qquad\qquad\qquad\qquad\qquad\qquad\qquad\qquad\qquad\qquad\qquad\qquad\qquad\qquad\qquad\]
\begin{equation}
=\frac{1}{s-1}\sum_{m=0}^{\infty}\frac{1}{(m+1)^k}\sum_{j=0}^{m+1}\binom{m+1}{j}(-1)^{j+1}\frac{1}{(x + j \ln a + (j + 1) \ln b)^{s-1}}.
\end{equation}
\begin{proof}
By using Lemma \ref{lem3}, we get
\begin{align*}
&=\frac{1}{\Gamma(s)}\int_0^{\infty}\frac{{\rm Li}_k(1 - e^{-t \ln ab})}{1 - e^{-t \ln ab}}e^{-t(x+\ln b)}t^{s-1}\int_{0}^{\ln ab}e^{-wt}dwdt\\
&=\frac{1}{\Gamma(s)}\int_0^{\infty}{\rm Li}_k(1 - e^{-t \ln ab})e^{-t(x+\ln b)}t^{s-2}dt\\
&=\frac{\Gamma(s-1)}{\Gamma(s)}(\xi_k (s - 1, x; a, b) - \xi_k (s - 1, x + \ln ab; a, b))\\
&=\frac{1}{s-1}\sum_{m=0}^{\infty}\frac{1}{(m+1)^k}\sum_{j=0}^{m+1}\binom{m+1}{j}(-1)^{j+1}\frac{1}{(x + j \ln a + (j + 1) \ln b)^{s-1}}.
\end{align*}
\end{proof}
So, we obtain the desired result.
\end{thm}

As a direct result of Raabe's formula and interpolation formula, we obtain following corollary for generalization of poly-Bernoulli numbers with $a, b$ parameters.

\begin{cor}
Raabe's formula in terms of generalization of poly-Bernoulli polynomials with
$a, b$ parameters is as follows:
\[\int_0^{\ln ab}B^{(k)}_n (-x - w; a, b) dw \qquad\qquad\qquad\qquad\qquad\qquad\qquad\qquad\qquad\qquad\qquad\qquad\qquad\qquad\qquad\qquad\qquad\qquad\qquad\qquad\qquad\qquad\]
\[=\frac{(-1)^{n+1}}{n+1}\sum_{m=0}^{\infty}\frac{1}{(m+1)^k}\sum_{j=0}^{m+1}\binom{m+1}{j}(-1)^{j+1}(x + j \ln a + (j + 1) \ln b)^{n+1}.\]
\end{cor}

\bigskip
\begin{flushleft}
{\bf Hassan Jolany}\\
Universit\'e des Sciences et Technologies de Lille\\
UFR de Math\'ematiques\\
Laboratoire Paul Painlev\'e\\
CNRS-UMR 8524 59655 Villeneuve d'Ascq Cedex/France\\
e-mail: hassan.jolany@math.univ-lille1.fr
\end{flushleft}

\bigskip
\begin{flushleft}
{\bf Roberto B. Corcino}\\
Mathematics and ICT Department\\
Cebu Normal University\\
Osmena Blvd., Cebu City\\
Philippines 6000\\
e-mail: rcorcino@yahoo.com
\end{flushleft}

\end{document}